\documentclass[12pt]{article}

\usepackage{amsfonts}

\newtheorem{thm}{Theorem}[section]

\newtheorem{cor}[thm]{Corollary}

\newcommand{\be}{\begin{equation}}
\newcommand{\ee}{\end{equation}}
\newcommand{\openbox}{\leavevmode
  \hbox to8pt{\hfil\vrule\vbox to6pt{\hrule width6pt\vfil\hrule}\vrule}}

\newcommand{\qed}{\hbox to5pt{ } \hfill \openbox\bigskip\medskip}

\newcommand{\cV}{\mbox{$\cal V$}}
\newcommand{\cC}{\mbox{$\cal C$}}
\newcommand{\cF}{\mbox{$\cal F$}}

\newcommand{\ve}[1]{\mathbf{#1}}

\title{Shattering bounds for tuple systems}
\author{G\'abor Heged\H{u}s
\\{\normalsize  \'Obuda University}
}
\begin{document}
\maketitle
\begin{abstract}
Let $V(n,d,q)$ stand for the $q$--ary Hamming spheres.
 Let  $\cC\subseteq (q)^n$ denote a tuple system such that 
$\cC\subseteq \cup_{i=0}^s \cV(n,d_i,q)$, where $d_1<\ldots <d_s$. 
We give here a general upper bound on the size of a shattered sets of the tuple system $\cC$.

\end{abstract}

\medskip
\section{Introduction}

First we introduce some notations.

Let  $[n]$ stand for the set $\{1,2,
\ldots, n\}$. We denote the family of all subsets of $[n]$  by $2^{[n]}$. 

Let $q\geq 2$ denote a fixed  integer throughout this paper. 
Let  $(q)$ stand for the set $\{0,1,
\ldots, q-1\}$.

Let $n>0$, $\cF\subseteq 2^{[n]}$ be a family of subsets of $[n]$, 
and $S$ be a subset of $[n]$. We say that $\cF$ {\it shatters} $S$ if
\begin{equation} \label{sha1}
\{F\cap S :~F\in \cF\}=2^{S}.
\end{equation}
Define
\begin{equation}
\mbox{sh} (\cF)=\{S\subseteq [n] :~\cF\mbox{\ shatters\ }S\}.
\end{equation}

The following well--known lemma was proved by 
Sauer in \cite{S}, and later independently by  Vapnik and Chervonenkis in \cite{VC},
and Perles and Shelah  in \cite{Sh}.

\begin{thm} \label{VC}
Suppose that $0\leq s \leq n-1$ and let $\cF\subseteq 2^{[n]}$ be 
an arbitrary set family with no shattered set of size $s+1$. Then
$$
|\cF|\leq \sum_{i=0}^{s} {n \choose i}.
$$
\end{thm}

Subsets  $\cV\subseteq (q)^n$ will be called {\em tuple
systems}. An element $\ve v$ of a tuple system $\cV$ can also be considered naturally as a 
function
from $[n]$ to $(q)$. 

Next we extends the usual 
binary notion of shattering. We say  
that the tuple  system $\cV$ {\em shatters} a set $S\subseteq [n]$, if 
$$
\{ \ve v\mid_{S}:~ \ve v\in \cV\}
$$
is the set of all functions from $S$ to $(q)$. Here $\ve v\mid_{S}$ denotes the 
restriction of the function  $\ve v$ to the set $S$. 

We define 
\begin{equation}
\mbox{Sh} (\cV):=\{S\subseteq [n]:~ \cV\mbox{\ shatters\ }S\},
\end{equation}
the set of the shattered sets of the tuple system $\cV$. 

Clearly $\mbox{Sh}(\cV)\subseteq 2^{[n]}$.
Moreover, if $\cF\subseteq 2^{[n]}$ is a set system, then we get 
$$
\mbox{sh}(\cF)= \mbox{Sh}(\{\ve v_F\in 2^{(n)}:~  F\in \cF\}).
$$

We can  also define the {\em $q$-ary Hamming spheres} $V(n,d,q)$.

Let $0\leq d\leq n$. Then
$$
\cV(n,d,q):=\{\ve v=(v_1,\ldots ,v_n)\in (q)^{[n]}:~ |\{i\in[n]:~ v_i\ne 0\}|=d\}.
$$

Our main result is similar to the Singleton Bound \cite{S2}.
\begin{thm} \label{main}
Suppose that the tuple system $\cC\subseteq (q)^n$ such that 
\begin{equation} \label{tart}
\cC\subseteq \cup_{i=0}^s \cV(n,d_i,q),
\end{equation}
 where $q\leq d_1<\ldots <d_s$. Define
 $t:=\max_i(d_{i+1}-d_i)$, and $k:=\min(d_s,n-t+1)$. Then
$|F|\leq k$ for each $F\in \mbox{sh}(\cC)$.
\end{thm}

Clearly our main result is sharp, namely we can choose $\cV:=\cup_{i=0}^s \cV(n,i,q)$, where $1\leq s\leq n$

Now we state our new upper bound for tuple systems.
\begin{cor} \label{main2}
Suppose that the tuple system $\cC\subseteq (q)^n$ such that 
\begin{equation} \label{tart2}
\cC\subseteq \cup_{i=0}^s \cV(n,d_i,q),
\end{equation}
 where $d_1<\ldots <d_s$. Define
 $t:=\max_i(d_{i+1}-d_i)$, and $k:=\min(d_s,n-t+1)$. Then
$$
|\cC|\leq \sum_{i=0}^k (q-1)^{n-i} {n \choose i}.
$$
\end{cor}

The choice $\cV:=\cup_{i=0}^s \cV(n,i,q)$, where $1\leq s\leq n$,
shows that this bound is sharp, too.

\section{Proofs}

We present here the proofs of our main results. 

{\bf Proof of Theorem \ref{main}:}
\smallskip

Let $m:=\mbox{max}\{|F|:~ F\in \mbox{Sh}(\cC)\}$. It is clear that $m\leq d_s$. Namely, proving in an indirect way, suppose that 
$m>d_s$. Then there exists an $F \in \mbox{Sh}(\cC)$
such that $|F|=m$, and there exists a vector $\ve v=(v_1,\ldots ,v_n)\in \cC$ with
$$
v_i=q-1\mbox{ for each }i\in F.
$$

Thus $|\{i\in [n]:~ v_i\neq 0\}|\geq |F|>d$, a contradiction with (\ref{tart}).  

Let $\mbox{supp}(\ve v):=\{i\in [n]:~ v_i\ne 0\}$ denote the support of the vector $\ve v$.

Now we can suppose that the distance $t$ appears between the levels 
$d_i$ and $d_{i+1}$, i.e., $d_{i+1}-d_i=t$.

First suppose that $d_i\geq m$. Then $d_i+t=d_{i+1}\leq n$, thus $n-t\geq d_i$.

Proving in an indirect way,  we suppose that $m\geq n-t+2$, then $n-t\leq m-2$, consequently 
$$
m\leq d_i\leq n-t\leq m-2,
$$
a contradiction.

Now consider the case $d_i+1\leq m$. Proving in an indirect way again,  we suppose that $n-t+2\leq m$. Clearly there exists an $F‌\in \mbox{Sh}(\cC)$ with $|F|=m$. Since $d_i+1\leq m$, hence there exists a $G\subseteq F$ with $|G|=d_i+1$. The condition $F\in \mbox{Sh}(\cC)$ means that there exists a vector $\ve v=(v_1,\ldots ,v_n)\in \cC$ with
\[
v_i=\left\{ \begin{array}{ll}
q-1 & \textrm{if $i\in G$} \\
0 & \textrm{if $i\in F\setminus G$.} 
\end{array} \right.
\]

Clearly 
\begin{equation} \label{szam}
|\mbox{supp}(\ve v)|\geq d_{i+1},
\end{equation}
because $G\subseteq \mbox{supp}(\ve v)$, $d_i+1=|G|\leq |\mbox{supp}(\ve v)|$, and  $\ve v\in \cC\subseteq \cup_{i=0}^s \cV(n,d_i,q)$. 

But then
$\mbox{supp}(\ve v)\subseteq G\cup ([n]\setminus F)$, hence 
$$
|\mbox{supp}(\ve v)|\leq |G\cup ([n]\setminus F)|=|G|+n-|F|=d_i+1+n-m
$$
and it follows from the indirect condition  that
$$
|\mbox{supp}(\ve v)|\leq d_i+1+t-2=d_i+t-1<d_{i+1},
$$
a contradiction with equation (\ref{szam}). \qed

Karpovsky and Milman  proved in \cite[Theorem 2]{KM} a natural generalization of 
Sauer's lemma for tuple systems. 

\begin{thm} \label{KaMi}
Let $0\leq s\leq n-1$ be an integer and let $\cV\subseteq (q)^n$ be a
tuple system with no shattered set of size $s+1$. 
Then
$$
|\cV|\leq \sum_{i=0}^s (q-1)^{n-i} {n \choose i}.
$$ 
\end{thm} \qed

{\bf Proof of Theorem \ref{main2}:}
\smallskip

It follows obviously from Theorem \ref{main} and \ref{KaMi}.

\end{document}